\documentclass[12pt]{article}

\usepackage{amsmath,amssymb,amsbsy,amsfonts,amsthm,latexsym,
                       amsopn,amstext,amsxtra,euscript,amscd}

\begin{document}
\newfont{\teneufm}{eufm10}
\newfont{\seveneufm}{eufm7}
\newfont{\fiveeufm}{eufm5}

\newfam\eufmfam
  \textfont\eufmfam=\teneufm \scriptfont\eufmfam=\seveneufm
  \scriptscriptfont\eufmfam=\fiveeufm
\def\frak#1{{\fam\eufmfam\relax#1}}

\def\bbbr{{\rm I\!R}} 
\def\bbbm{{\rm I\!M}}
\def\bbbn{{\rm I\!N}} 
\def\bbbf{{\rm I\!F}}
\def\bbbh{{\rm I\!H}}
\def\bbbk{{\rm I\!K}}
\def\bbbp{{\rm I\!P}}
\def\bbbone{{\mathchoice {\rm 1\mskip-4mu l} {\rm 1\mskip-4mu l}
{\rm 1\mskip-4.5mu l} {\rm 1\mskip-5mu l}}}
\def\bbbc{{\mathchoice {\setbox0=\hbox{$\displaystyle\rm C$}\hbox{\hbox
to0pt{\kern0.4\wd0\vrule height0.9\ht0\hss}\box0}}
{\setbox0=\hbox{$\textstyle\rm C$}\hbox{\hbox
to0pt{\kern0.4\wd0\vrule height0.9\ht0\hss}\box0}}
{\setbox0=\hbox{$\scriptstyle\rm C$}\hbox{\hbox
to0pt{\kern0.4\wd0\vrule height0.9\ht0\hss}\box0}}
{\setbox0=\hbox{$\scriptscriptstyle\rm C$}\hbox{\hbox
to0pt{\kern0.4\wd0\vrule height0.9\ht0\hss}\box0}}}}
\def\bbbq{{\mathchoice {\setbox0=\hbox{$\displaystyle\rm
Q$}\hbox{\raise
0.15\ht0\hbox to0pt{\kern0.4\wd0\vrule height0.8\ht0\hss}\box0}}
{\setbox0=\hbox{$\textstyle\rm Q$}\hbox{\raise
0.15\ht0\hbox to0pt{\kern0.4\wd0\vrule height0.8\ht0\hss}\box0}}
{\setbox0=\hbox{$\scriptstyle\rm Q$}\hbox{\raise
0.15\ht0\hbox to0pt{\kern0.4\wd0\vrule height0.7\ht0\hss}\box0}}
{\setbox0=\hbox{$\scriptscriptstyle\rm Q$}\hbox{\raise
0.15\ht0\hbox to0pt{\kern0.4\wd0\vrule height0.7\ht0\hss}\box0}}}}
\def\bbbt{{\mathchoice {\setbox0=\hbox{$\displaystyle\rm
T$}\hbox{\hbox to0pt{\kern0.3\wd0\vrule height0.9\ht0\hss}\box0}}
{\setbox0=\hbox{$\textstyle\rm T$}\hbox{\hbox
to0pt{\kern0.3\wd0\vrule height0.9\ht0\hss}\box0}}
{\setbox0=\hbox{$\scriptstyle\rm T$}\hbox{\hbox
to0pt{\kern0.3\wd0\vrule height0.9\ht0\hss}\box0}}
{\setbox0=\hbox{$\scriptscriptstyle\rm T$}\hbox{\hbox
to0pt{\kern0.3\wd0\vrule height0.9\ht0\hss}\box0}}}}
\def\bbbs{{\mathchoice
{\setbox0=\hbox{$\displaystyle     \rm S$}\hbox{\raise0.5\ht0\hbox
to0pt{\kern0.35\wd0\vrule height0.45\ht0\hss}\hbox
to0pt{\kern0.55\wd0\vrule height0.5\ht0\hss}\box0}}
{\setbox0=\hbox{$\textstyle        \rm S$}\hbox{\raise0.5\ht0\hbox
to0pt{\kern0.35\wd0\vrule height0.45\ht0\hss}\hbox
to0pt{\kern0.55\wd0\vrule height0.5\ht0\hss}\box0}}
{\setbox0=\hbox{$\scriptstyle      \rm S$}\hbox{\raise0.5\ht0\hbox
to0pt{\kern0.35\wd0\vrule height0.45\ht0\hss}\raise0.05\ht0\hbox
to0pt{\kern0.5\wd0\vrule height0.45\ht0\hss}\box0}}
{\setbox0=\hbox{$\scriptscriptstyle\rm S$}\hbox{\raise0.5\ht0\hbox
to0pt{\kern0.4\wd0\vrule height0.45\ht0\hss}\raise0.05\ht0\hbox
to0pt{\kern0.55\wd0\vrule height0.45\ht0\hss}\box0}}}}
\def\bbbz{{\mathchoice {\hbox{$\sf\textstyle Z\kern-0.4em Z$}}
{\hbox{$\sf\textstyle Z\kern-0.4em Z$}}
{\hbox{$\sf\scriptstyle Z\kern-0.3em Z$}}
{\hbox{$\sf\scriptscriptstyle Z\kern-0.2em Z$}}}}
\def\ts{\thinspace}

\newtheorem{theorem}{Theorem}[section]
\newtheorem{corollary}{Corollary}
\newtheorem*{main}{Main Theorem}
\newtheorem{lemma}[theorem]{Lemma}
\newtheorem{proposition}{Proposition}
\newtheorem{conjecture}{Conjecture}
\newtheorem*{problem}{Problem}
\theoremstyle{definition}
\newtheorem{definition}[theorem]{Definition}
\newtheorem{remark}{Remark}
\newtheorem*{notation}{Notation}
\newcommand{\ep}{\varepsilon}
\newcommand{\eps}[1]{{#1}_{\varepsilon}}

\def\squareforqed{\hbox{\rlap{$\sqcap$}$\sqcup$}}
\def\qed{\ifmmode\squareforqed\else{\unskip\nobreak\hfil
\penalty50\hskip1em\null\nobreak\hfil\squareforqed
\parfillskip=0pt\finalhyphendemerits=0\endgraf}\fi}

\def\eqref#1{(\ref{#1})}

\def\cA{{\mathcal A}}
\def\cB{{\mathcal B}}
\def\cC{{\mathcal C}}
\def\cD{{\mathcal D}}
\def\cE{{\mathcal E}}
\def\cF{{\mathcal F}}
\def\cG{{\mathcal G}}
\def\cH{{\mathcal H}}
\def\cI{{\mathcal I}}
\def\cJ{{\mathcal J}}
\def\cK{{\mathcal K}}
\def\cL{{\mathcal L}}
\def\cM{{\mathcal M}}
\def\cN{{\mathcal N}}
\def\cO{{\mathcal O}}
\def\cP{{\mathcal P}}
\def\cQ{{\mathcal Q}}
\def\cR{{\mathcal R}}
\def\cS{{\mathcal S}}
\def\cT{{\mathcal T}}
\def\cU{{\mathcal U}}
\def\cV{{\mathcal V}}
\def\cW{{\mathcal W}}
\def\cX{{\mathcal X}}
\def\cY{{\mathcal Y}}
\def\cZ{{\mathcal Z}}




\newcommand{\comm}[1]{\marginpar{\fbox{#1}}}

\newcommand{\ignore}[1]{}

\def\vec#1{\mathbf{#1}}

\def\e{\mathbf{e}}



\def\AA{\mathbb{A}}
\def\BB{\mathbf{B}}

\def\dist{\mathrm{dist}}

\hyphenation{re-pub-lished}

\def\vol{{\mathrm{vol}\,}}
\def\ad{{\mathrm ad}}

\def \F{{\bbbf}}
\def \K{{\bbbk}}
\def \nd{{\, | \hspace{-1.5 mm}/\,}}

\def \Z{{\bbbz}}
\def\Zn{\Z_n}
\def \N{{\bbbn}}
\def \Q{{\bbbq}}
\def \R{{\bbbr}}
\def\Fp{\F_p}
\def \fp{\Fp^*}
\def\\{\cr}
\def\({\left(}
\def\){\right)}
\def\fl#1{\left\lfloor#1\right\rfloor}
\def\rf#1{\left\lceil#1\right\rceil}

\def\Spec#1{\mbox{\rm {Spec}}\,#1}
\def\invp#1{\mbox{\rm {inv}}_p\,#1}
\def\ADM{\'{A}d\'{a}m}

\def\ADMPR {\ADM\  property}
\def\SADMPR {spectral \ADM\  property}
\def\CG{circulant  graph}
\def\AM{adjacency matrix}
\def\AMs{adjacency matrices}

\def\Ln#1{\mbox{\rm {Ln}}\,#1}

\def\epp{\mbox{\bf{e}}_{p-1}}
\def\ep{\mbox{\bf{e}}_p}
\def\ed{\mbox{\bf{e}}_{d}}

\def\ii {\iota}

\def\wt#1{\mbox{\rm {wt}}\,#1}

\def\GR#1{{ \langle #1 \rangle_n }}

\def\ab{\{\pm a,\pm b\}}
\def\cd{\{\pm c,\pm d\}}

\def\Bt {\mbox{\rm {Bt}}}

\def\Res#1{\mbox{\rm {Res}}\,#1}

\def\Tr#1{\mbox{\rm {Tr}}\,#1}


\title{On Some Dynamical
Systems in  Finite Fields\\ and Residue Rings}

\author{ 
{\sc Igor E. Shparlinski} \\
{Department of Computing, Macquarie University} \\
{Sydney, NSW 2109, Australia} \\
{igor@ics.mq.edu.au}}

\maketitle

\begin{abstract}
We use character sums to confirm  several recent conjectures of
V.~I.~Arnold  on the uniformity of distribution properties of
a certain dynamical system in a finite field. On the other hand,
we show that some conjectures are wrong. We also analyze several
other conjectures of
V.~I.~Arnold  related to the orbit length of similar dynamical
systems in residue rings and outline possible ways to prove
them. We also show that some of them require
further tuning.
\end{abstract}

\section{Introduction}

In a recent series of papers,
V.~I.~Arnold~\cite{Arn1,Arn2,Arn3,Arn4,Arn5,Arn6}
has considered dynamical systems related to linear transformations
in finite fields and residue rings and made a number of
conjectures. We observe
that  the study of the length, distribution of element and other
properties, of the orbits of such dynamical systems
has a long and successful history, which dates back to
early works of N.~M.~Korobov~\cite{Kor},
H.~Niederreiter~\cite{Nied1,Nied5},  A.~G.~Postnikov~\cite{Post}
and many other researchers. Here we show that some classical results
immediately imply some of these conjectures. We  also show that
several other conjectures are not correct as they are stated
in~\cite{Arn1,Arn2,Arn3,Arn4,Arn5,Arn6} and need some adjustments.

For a prime $p$ and a positive integer $n$,  we
denote by  $\F_{p^n}$ the finite field of $p^n$
elements (we refer to~\cite{LN} for the background
information on finite fields).

  We   fix a primitive root
$\vartheta$ of
$\F_{p^n}$ and recall that we also have $\F_{p^n} \cong \F_p[\vartheta]$.
In particular, $\F_{p^n}$ can be considered as an $n$-dimensional
vector space over $\F_p$, where with each element
$\alpha \in \F_{p^n}$ one can associate the coordinate vector
$\vec{a} = (a_0, \ldots, a_{n-1}) \in \F_p^n$ from the expansion
$$
\alpha = \sum_{j=0}^{n-1} a_j \vartheta^{j-1}.
$$
Accordingly, V.~I.~Arnold~\cite{Arn4} suggests
to study the sequence of vectors $\vec{a}_m = (a_{0,m}, \ldots, a_{n-1,m}) \in
\F_p^n$  corresponding to the powers
\begin{equation}
\label{eq:Vector}
\vartheta^m = \sum_{j=0}^{n-1} a_{j,m} \vartheta^{j-1}.
\end{equation}

Clearly, assuming that   $\F_p$ is represented
by the elements of the set $\{0,1, \ldots, p-1\}$,
one can view the points
\begin{equation}
\label{eq:Points}
\frac{1}{p} \vec{a}_m, \qquad m =1,\ldots, M,
\end{equation}
as $M$ points of an $n$-dimensional unit cube $[0,1]^n$.
For $M = p^n-1$ these points form a regular cubic lattice
(with only one missing point $(0, \ldots, 0)$).
It has also been conjectured by V.~I.~Arnold~\cite{Arn4}
that in fact even the first $M <  p^n-1$ powers
already form a rather uniformly distributed point
set.
  Namely, given a region $\Omega \in  [0,1]^n$ with  smooth
boundary, we denote by $N_\vartheta(M,\Omega)$ the number of
points~\eqref{eq:Points} which belong to $\Omega$.
The conjecture of Section~2.A of~\cite{Arn4} asserts that
\begin{equation}
\label{eq:Conj 2A}
N_\vartheta(M,\Omega) = M \vol \Omega + o(M)
\end{equation}
provided that   $M \sim \mu p^n$ for
some fixed $\mu>0$ (and $p \to \infty$).

We start with an observation that using classical bounds of incomplete
exponential sums with exponential functions, see~\cite{KoSh,Kor,LN},
and some standard tools
from  the theory of uniform distribution,
  see  Section~\ref{sec:discr},  one can
derive the following
improved version of the conjecture~\eqref{eq:Conj 2A}:
\begin{equation}
\label{eq:Conj 2A-strong}
N_\vartheta(M,\Omega) = M \vol \Omega +  O\(M^{1-1/n} p^{1/2n} (\log
p)^{1+1/n}\),
\end{equation}
which is nontrivial whenever $M/p^{1/2} (\log p)^{n+1} \to \infty$.
In fact, using some results of H. Niederreiter~\cite{Nied1,Nied5}
one can easily extend the above result in several directions.

In fact, using the
results of  J.~Bourgain and M.-C.~Chang~\cite{BourChang},
which in turn generalize
several recently emerged results of J.~Bourgain,  A.~A.~Glibichuk and
S.~V.~Konyagin~\cite{BGK,BourKon},
one can also study the distribution in intervals
of the set~\eqref{eq:Points}  for
extremely small values of $M$.
For example, see~\cite{Chang} for more details and
a version of the bound~\eqref{eq:Conj 2A-strong} which is nontrivial
provided that $M \ge p^\varepsilon$ for any fixed $\varepsilon > 0$
and sufficiently large $p$. The bound of the error term in~\cite{Chang}
is not completely explicit, so for large values of $M$ the
bound~\eqref{eq:Conj 2A-strong}  is better than that of~\cite{Chang}.

Moreover, for $n=1$, that is, for prime fields, using bounds of
exponential sums from~\cite{BGLS,Gar2,GarShp}, one can obtain nontrivial
results for even smaller intervals, which however holds only
for almost all primes $p$ (rather than for all $p$).

Furthermore, motivated by the results of~\cite{Bour1,FuWe},
we consider the distribution of vectors $\vec{a}_m$ where instead of
an initial segment $[1,M]$,  $m$ runs
through the values of a polynomial. Unfortunately, we are not able to
treat arbitrary polynomials with integer coefficients for every primitive root
$\vartheta$  but rather obtain a result which holds for almost all
primitive roots. However,  in the case of monomials,
employing   the bound of exponential sums with
the sequence $\vartheta^{m^k}$, $m =1, 2, \ldots$,  from~\cite{FHS1},
we obtain a nontrivial estimate for every  primitive root
$\vartheta$.

V.~I.~Arnold~\cite{Arn1,Arn3,Arn5,Arn6} also describes similar
dynamical systems in the residue ring $\Z_\ell$ modulo $\ell$ and makes several
conjectures  about the length of the orbits. More specifically,
given an integer $g\ge 2$ with $\gcd(g, \ell)=1$,
V.~I.~Arnold~\cite{Arn1,Arn3,Arn5,Arn6} suggests to consider
the dynamical properties of the residues $g^m \pmod \ell$.

We recall that the {\it Carmichael function\/} $\lambda(\ell)$
is defined for all $\ell \ge 1$ as the largest order of any
element in the multiplicative group $\Z_\ell^*$.
More explicitly, for any prime power $p^\nu$, one has
$$
\lambda(p^\nu) = \left\{  \begin{array}{ll}
p^{\nu-1}(p-1)&\qquad\hbox{if $p \ge 3$ or $\nu \le 2$},\\
2^{\nu-2}&\qquad\hbox{if $p = 2$ and $\nu \ge 3$},
\end{array} \right.
$$
and for an arbitrary integer $\ell\ge 2$,
$$
\lambda(\ell)=\mathrm{lcm}\(\lambda( p_1^{\nu_1}),\ldots,
\lambda(p_s^{\nu_s})\),
$$
where $\ell=p_1^{\nu_1}\ldots p_s^{\nu_s}$ is the prime factorization
of $\ell$. Clearly, $\lambda(1)=1$.
We also let $\varphi(\ell)$ denote the {\it Euler function\/},
which is defined as usual by
$$
\varphi(\ell)=\#\Z_\ell^* =\prod_{j=1}^sp_j^{\nu-1}(p_j-1),
$$
with $\varphi(1) =1$.
  Finally, for a an integer $g$ with $\gcd(g,
\ell)=1$, we denote by $t_g(\ell)$ the multiplicative order of $g$
modulo $\ell$.
Clearly, we have the divisibilities
$$
t_g(\ell)~|~\lambda(\ell)~|~\varphi(\ell).
$$
Several conjectures of~\cite{Arn1,Arn3,Arn5,Arn6}
can be reformulated as various statements about the relative size
of $t_g(\ell)$, $\lambda(\ell)$ and $\varphi(\ell)$, on average and
individually. We discuss these conjectures and show that
some of them are already known in the literature, while some
can be proved to be  wrong.
It is suggested in Section~1 of~\cite{Arn5}
that for $g =2$ the average multiplicative order
$$
T_g(L) =
\frac{1}{L} \sum_{\substack{\ell =1\\ \gcd(g,\ell)=1}}^L t_g(\ell)
$$
grows like
\begin{equation}
\label{eq:Conj Av Period}
T_g(L) \sim c(g) \frac{L}{\log L}
\end{equation}
for some constant $c(g)>0$ depending only on $g$ (we note
that in~\cite{Arn5} it is made explicit only for $g = 2$).

We show that   the classical result of Hooley~\cite{Hool}
on {\it Artin's conjecture\/}, implies,  under the
{\it Extended Riemann Hypothesis\/}, that
the conjecture~\eqref{eq:Conj Av Period} is wrong and in fact
$$
T_g(L) \ge  \frac{L}{\log L}
\exp\(C(g)(\log \log \log L)^{3/2} \).
$$
for some constant $C(g)>0$ depending only on $g$.
Furthermore, we believe that in fact $T_g(L)$ grows
even faster. It is possible that  the
method of proof of  Theorems~1 and~2 in~\cite{BFLPS}, which in turn
is an extension
of the method of~\cite{MaPom} (see also~\cite{Ford1}),
together with the result of Hooley~\cite{Hool},  can be used to derive
that, under the  Extended Riemann Hypothesis,
\begin{equation}
\label{eq:AverPerTight}
T_g(L) \ge  \frac{L}{\log L}  \exp\( (\log \log \log L)^{2+o(1)}\).
\end{equation}
For the upper bound it is probably natural to assume that
\begin{equation}
\label{eq:T and lambda}
T_g(L) = o\(\frac{1}{L}  \sum_{\ell =1}^L \lambda (\ell)\) .
\end{equation}
Note the sum on the right hand side of~\eqref{eq:T and lambda}
has been estimated by P.~Erd{\H o}s, C.~Pomerance and
E.~Schmutz~\cite{EPS}.

Finally, we give a guide to the
literature concerning results and methods which can probably be
of great use for the theory of algebraic dynamical
systems over finite fields and rings.

It is very well known that there are close ties between
  number-theory  and dynamical systems. For
example, one can associate dynamical systems with continued fractions,
various number systems, the $3x+1$ transformation and other
number-theoretic constructions.
A wealth of very interesting results can be found in the literature.

However, we would like to use this paper as an opportunity to attract
more attention
of the dynamical system community to a great variety of
already existing
number theoretic results and techniques which can be
of great significance for  studying various algebraic
dynamical systems.  In
particular, these include, but are not limited too, bounds on
various exponential sums, periods of
various sequences and average  values of associated arithmetic functions,
For this very purpose we
do not try to formulate and prove our results in their full generality
but rather limit ourselves to the most interesting and illuminating special
cases.   We however indicate possible extensions of our results
and directions for further research.

Although exponential
sums have been used for this purpose, see the work of
M. Degli Esposti and S. Isola~\cite{DEI} and of P.~Kurlberg
and Z.~Rudnick~\cite{KuRu4},
their full potential seems to be not fully used in the dynamical system
theory.  We would like to stress that all such applications follow the
same pattern:
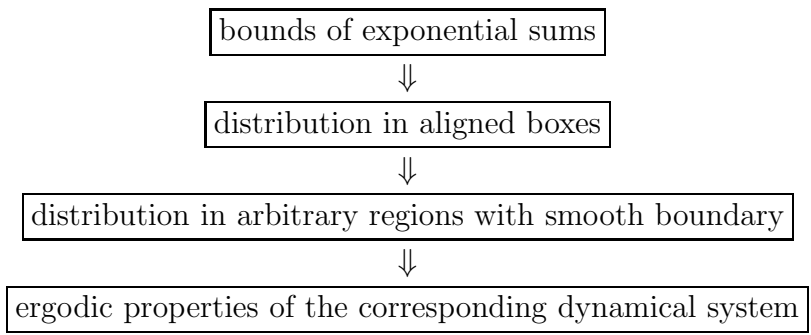
\begin{eqnarray*}
\centerline{\fbox{bounds of exponential sums}}\\
\centerline{$\Downarrow$}\\
\centerline{\fbox{distribution in aligned boxes}}\\
\centerline{$\Downarrow$}\\
\centerline{\fbox{distribution in arbitrary regions with  smooth
boundary}}\\
\centerline{$\Downarrow$}\\
\centerline{\fbox{ergodic properties of the corresponding dynamical system}}
\end{eqnarray*}

The link between exponential sums and the distribution in aligned boxes
is provided by the Koksma--Sz\"usz inequality, see  Theorem~1.21
of~\cite{DrTi}.

The link between the distribution in aligned boxes and arbitrary
regions is given by the results of H.~Niederreiter and J.~M.~Wills~\cite{NiWi}
and their more recent refinement of
M.~Laczkovich~\cite{Lac}.

Surprisingly enough, the essentially tautological link between the
distribution in  arbitrary  regions and ergodic properties has never been
exploited in a  systematic way, although it definitely deserves much
more attention
which we hope to attract with this paper.

\section{Dynamical systems in finite fields}

\subsection{Preliminaries}

Here we show how well known bounds of exponential sums
can be used to derive various results about the orbits of
$\vartheta^{f(m)}$, with a polynomial $f(X) \in \Z[X]$.

Throughout this section, any implied constants
in the symbols $O$ may depend on $n$ and $\Omega$ (and occasionally, where
obvious, on an integer parameter $k$).

As we have mentioned, the results of in this section can be extended in
several  directions.

\subsection{Background on finite fields}

Let $\omega_0, \ldots, \omega_{n-1}$ be a basis of $\F_{p^n}$ over $\F_p$
which is  dual to the basis $1, \vartheta, \ldots, \vartheta^{n-1}$.
That is,
$$
\Tr\(\omega_i \vartheta^j\)= \left\{  \begin{array}{ll}
1 ,& \quad \mbox{if}\  i=j;\\
0, & \quad \mbox{otherwise};
\end{array} \right.
\qquad 0 \le i,j\le  n-1,
$$
where
$$
\Tr(\alpha) = \sum_{k =0}^{n-1} \alpha^{p^k}
$$
is the trace of $\alpha \in   \F_{p^n}$ in $\F_p$.

Therefore, from~\eqref{eq:Vector} we derive
\begin{equation}
\label{eq:Basis Rep}
a_{j,m} = \Tr\(\omega_j \vartheta^m\)
\end{equation}
for every $j =0, \ldots n-1$ and $m =1,2, \ldots$.

It is also useful to recall that there are $\varphi(p^n-1)$ primitive roots of
$\F_{p^n}$.

\subsection{Background on exponential sums}

Let us denote $\ep(z)  = \exp(2 \pi i z/p)$. Then for every $\gamma
\in \F_{p^n}$
the function
$\alpha \mapsto \ep\(\gamma \Tr(\alpha)\)$ is  an additive character
of $\F_{p^n}$.

  For example,
it follows immediately from a combination of Theorem~8.24 and
Theorem~8.81 of~\cite{LN} (see also~\cite{KoSh,Kor} and the references
therein), that for any $M \le p^n-2$
the following bound holds
\begin{equation}
\label{eq:Classic Exp Sum}
\max_{\gamma \in \F_{p^n}^*}
  \left| \sum_{m=1}^M \ep\(\Tr(\gamma \vartheta^m)\)\right| = O(p^{n/2} \log p).
\end{equation}

The following estimate is a special case of a more general result
of~\cite{GarKar}.
We remark that in~\cite{GarKar} it is shown only in the case $n=1$ but the
proof extends to arbitrary fields without any changes.

\begin{lemma}
\label{lem:Double Exp Sum} For any primitive root $\vartheta \in \F_{p^n}$,
any two subsets $\cX, \cY \in \Z_{p^n-1}$ and any function $\psi(y)$ with
$$
\max_{y \in \cY} |\psi(y)| \le \Psi,
$$
the following bound holds
$$
\max_{\gamma \in \F_{p^n}^*}
\sum_{x \in \cX}  \left| \sum_{y \in \cY} \psi(y) \ep\(\Tr(\gamma
\vartheta^{xy})\)\right| =
O\(\Psi p^{9n/8 +o(1)} \(\#\cY\)^{3/4}\).
$$
\end{lemma}

We know recall the bound of exponential sums with  $\vartheta^{m^k}$
from~\cite{FHS1} which we use in the proof of Theorem~\ref{thm:Ergod-k}.
More precisely, we use Theorem~6 (for $k=2$) and Theorem~7 (for $k\ge 3$)
of~\cite{FHS1} (we remark that in~\cite{FHS1} these results are proven
only for $n=1$ but the general case can be obtained by a simple
typographical change of $p$ to $p^n$).

Let us define
\begin{equation}
\label{eq:rho}
\rho(k) =  \left\{  \begin{array}{ll}
\displaystyle
\frac{1}{8} ,& \quad \mbox{if}\  k =2;\\
\\
\displaystyle
\frac{ \rf{k/2} -1 }{2k  \rf{k/2} + 2} , & \quad \mbox{if}\ k \ge
3.
\end{array} \right.
\end{equation}

\begin{lemma}
\label{lem:Exp Sum-k} For  any primitive root $\vartheta \in \F_{p^n}$,
the following bound holds
$$
\max_{\gamma \in \F_{p^n}^*}
  \left| \sum_{m=1}^{p^n-1} \ep\(\Tr(\gamma \vartheta^{m^k})\)\right| \le
p^{(1-\rho(k))n + o(1)} .
$$
\end{lemma}

\subsection{Background on discrepancies}
\label{sec:discr}

For a finite set $\cU \subseteq [0,1]^n$ and domain  $\Omega \subseteq
[0,1]^n $, we define the {\it $\Omega$-discrepancy\/}
$$
\Delta(\cU,\Omega) = \left| \frac{ \#\{ \vec{u} \in\cU \cap \Omega\}
}{\#\cU } - \vol \Omega\right|,
$$
and  the {\it box discrepancy\/} of $\cU$,
$$
D(\cU) = \sup_{\BB  \subseteq [0,1]^n}  \Delta(\cU,\BB) ,
$$
where the supremum is taken over all boxes $\BB = [\alpha_1,
\beta_1] \times \ldots \times [\alpha_k, \beta_k]$.

We define the  distance between a vector $\vec{u} \in [0,1]^n$ and a
set $\Gamma \subseteq [0,1]^n $  by
$$
  \dist(\vec{u},\Gamma) = \inf_{\vec{w} \in\Gamma}
\|\vec{u} - \vec{w}\|
$$
where $\|\vec{v}\|$ denotes the Euclidean norm of $\vec{v} \in \R^n$. Given
$\varepsilon >0$ and a domain  $\Omega \subseteq [0,1]^n $ we define
the  sets
$$
\Omega_\varepsilon^{+} = \left\{ \vec{u} \in  [0,1]^n \backslash
\Omega \ | \ \dist(\vec{u},\Omega) < \varepsilon \right\}
$$
and
$$
\Omega_\varepsilon^{-} = \left\{ \vec{u} \in \Omega \ | \
\dist(\vec{u},[0,1]^n \backslash \Omega )  < \varepsilon  \right\} .
$$

Let $b(\varepsilon)$ be any increasing function defined for
$\varepsilon >0$ and such that $\lim_{ \varepsilon \to
0}b(\varepsilon) = 0$. Following~\cite{Lac,NiWi}, we define the
class $\cM_b$ of  domains  $\Omega \subseteq [0,1]^n $ for which
$$
\vol \Omega_\varepsilon^{+}  \le b(\varepsilon)
\qquad \mbox{and}
\qquad
\vol \Omega_\varepsilon^{-}  \le b(\varepsilon) .
$$

As special case of a result of H.~Weyl~\cite{Weyl}
implies that for domains $\Omega$ with a piecewise smooth
boundary, one can take $b(\varepsilon) = O( \varepsilon)$.

\begin{lemma}
\label{lem:Omega-eps} For any domain  $\Omega \subseteq [0,1]^n$
with a piecewise smooth boundary
$$
\vol \Omega_\varepsilon^{\pm } = O\(\varepsilon\)
$$
\end{lemma}

A relation  between $D(\cU)$ and $\Delta(\cU,\Omega)$ for $\Omega \in
\cM_b$ is given by the following inequality from~\cite{Lac} (see
also~\cite{NiWi}).

\begin{lemma}
\label{lem:Omega Discr} For any domain  $\Omega \in \cM_b$, we have
$$
\Delta(\cU,\Omega)  = O \( b\(n^{1/2} D(\cU)^{1/n}\)\) .
$$
\end{lemma}

The {\it Koksma--Sz\"usz inequality\/}, see
Theorem~1.21 of~\cite{DrTi},   provides an important link
between box discrepancy and exponential sums:

\begin{lemma}
\label{le:K-S} For any integer $L > 1$, and a set $ \cU \subseteq [0,1]^n$
     of $M$ points, one has
$$
D(\cU)  = O\(  \frac{1 }{ L } + \frac{1}{M} \sum_{\substack{\vec{c} =
(c_0, \ldots, c_{n-1}) \in \Z^n\setminus \{ \vec{0}\} \\ |c_j| \le  L
\ j=0, \ldots, n-1}} \prod_{j=0}^{n-1} \frac{1}{ (1+|c_j|)} \left|
\sum_{\vec{u}\in \cU} \exp \( 2 \pi i \vec{c}\cdot\vec{u} \) \right| \) ,
$$
where
$$\vec{c}\cdot\vec{u} = \sum_{j=0}^{n-1} c_j u_j
$$
denotes the inner product of $\vec{c} = (c_0, \ldots, c_{n-1})$
and $\vec{u} = (u_0, \ldots, u_{n-1})$.
\end{lemma}

To estimate the box discrepancy of the set~\eqref{eq:Points}
we apply Lemma~\ref{le:K-S} with $L = (p-1)/2$.
By~\eqref{eq:Basis Rep} we see that the corresponding exponential
sums takes shape
$$
\sum_{m=1}^M \ep\(  \sum_{j=0}^{n-1} c_j \Tr(\omega_j \vartheta^m)\)
= \sum_{m=1}^M \ep\(   \Tr(\gamma \vartheta^m)\)
$$
where $\gamma = c_0 \omega_0 + \ldots + c_{n-1} \omega_{n-1} \in \F_{p^n}^*$.

  Applying the bound~\eqref{eq:Classic Exp Sum} together with
Lemma~\ref{le:K-S},
we see that the box discrepancy of the set~\eqref{eq:Points} is
$O\(M^{-1} p^{n/2} (\log p)^{n+1}\)$, see also~\cite{Nied1,Nied5}
and references therein
for several   more general results.
  Now the
bound~\eqref{eq:Conj 2A-strong} follows
directly from  Lemmas~\ref{lem:Omega-eps} and~\ref{lem:Omega Discr}.

\subsection{Distribution of points in orbits}.

For an polynomial $f(X) \in \Z[X]$
and a given region $\Omega \in  [0,1]^n$ with  smooth
boundary, we denote by $N_{\vartheta}(f;M,\Omega)$ the number of
points
\begin{equation}
\label{eq:Points-poly}
\frac{1}{p} \vec{a}_{f(m)}, \qquad m =1,\ldots, M,
\end{equation}
which belong to $\Omega$.

\begin{theorem}
\label{thm:Ergod-f} Let $f(X) \in \Z[X]$ be a fixed nonconstant
polynomial and let $\varTheta$
be the set of all
$\varphi(p^n-1)$ primitive roots of
$\F_{p^n}$.  For  any positive integer
$M
\le p^n-1$  and any
  region $\Omega \in [0,1]^n$ with piecewise  smooth
boundary, we have
$$
\frac{1}{\varphi(p^n-1)}\sum_{\vartheta \in \varTheta}
\left|N_{\vartheta}(f;M,\Omega)  - M
\vol
\Omega
\right|  \le M^{1-1/4n} p^{1/8 + o(1)}.
$$
\end{theorem}

\begin{proof}
To estimate the box discrepancy of the set~\eqref{eq:Points-poly}
we apply Lemma~\ref{le:K-S} with $L = (p-1)/2$.
As in Section~\ref{sec:discr}, by~\eqref{eq:Basis Rep} we see that
the corresponding
exponential sums takes shape
$$
\sum_{m=1}^M \ep\(  \sum_{j=0}^{n-1} c_j \Tr\(\omega_j \vartheta^{f(m)}\)\)
= \sum_{m=1}^M \ep\(    \Tr\(\gamma \vartheta^{f(m)}\)\)
$$
where $\gamma = c_0 \omega_0 + \ldots + c_{n-1} \omega_{n-1} \in \F_{p^n}^*$.
Applying   Lemma~\ref{le:K-S},
we see that the box discrepancy $D_\vartheta(f;M)$ of the set~\eqref{eq:Points}
satisfies
$$
\sum_{\vartheta \in \varTheta} D_\vartheta(f;M) = O\(\frac{1}{p} +
\frac{(\log p)^n}{M}
\max_{\gamma \in \F_{p^n}^*} \sum_{\vartheta \in \varTheta}
  \left| \sum_{m=1}^{M} \ep\(\Tr(\gamma \vartheta^{f(m)})\)\right| \).
$$
Let $\cX$ be the set of all elements of $\Z_{p^n-1}$ which are relatively
prime to $p^n-1$. Fix an arbitrary primitive root $\vartheta_0 \in \F_{p^n}$.
Then $\varTheta = \{ \vartheta_0^x\ | \ x \in \cX\}$.

We denote by $\cY$ the value set $\cY = \{f(m) \pmod {p^n-1}\ | \ m
=1, \ldots, M\}$
and by $\psi(y)$ the multiplicity of $y \in \cY$ (that is, the number of
$m =1, \ldots, M$
with $y \equiv f(m) \pmod {p^n-1}$). In particular,
$\# \cY \le M$ and by the  famous
Nagell--Ore theorem (see~\cite{Huxley} for its strongest known form) we have
$$
\Psi = \max_{y \in \cY} |\psi(y)| = p^{o(1)}.
$$
We derive from
Lemma~\ref{lem:Double Exp Sum} that
$$
\max_{\gamma \in \F_{p^n}^*} \sum_{\vartheta \in \varTheta}
  \left| \sum_{m=1}^{M} \ep\(\Tr(\gamma \vartheta^{f(m)})\)\right|
\le  p^{9n/8 +o(1)} M^{3/4},
$$
which implies the bound
\begin{equation}
\label{eq:PolyDiscr-Aver}
\sum_{\vartheta \in \varTheta} D_\vartheta(f;M) \le  p^{9n/8 +o(1)} M^{-1/4}.
\end{equation}
  Now,
since $\Omega$ has a piecewise smooth boundary,
from  Lemmas~\ref{lem:Omega-eps}
and~\ref{lem:Omega Discr} and the H{\"o}lder inequality,  we derive
  \begin{eqnarray*}
\lefteqn{\sum_{\vartheta \in \varTheta} \left|N_{\vartheta}(f;M,\Omega)  - M
\vol
\Omega
\right| = O\(M \sum_{\vartheta \in \varTheta} D_\vartheta(f;M)^{1/n}\)}\\
& &  \qquad \qquad\qquad = O\(M \(\varphi(p^n-1)^{n-1}\sum_{\vartheta
\in \varTheta}
D_\vartheta(f;M)\)^{1/n}\)\\
   & &  \qquad \qquad\qquad  = O\(M
\varphi(p^n-1)^{1-1/n}\(\sum_{\vartheta \in \varTheta}
D_\vartheta(f;M)\)^{1/n}\).
\end{eqnarray*}
Since $k/\varphi(k) = O(\log \log k)$ for every integer $k$, see Theorem~328
of~\cite{HardyWright}, from~\eqref{eq:PolyDiscr-Aver} we derive the
desired estimate.
\end{proof}

  For example, we see from Theorem~\ref{thm:Ergod-f} that for every fixed
$\varepsilon$ and $M \ge p^{n/2 + \varepsilon}$, for almost all
  primitive roots of
$\vartheta \in \F_{p^n}$, we have $N_{\vartheta}(f;M,\Omega)  \sim M
\vol \Omega$ for every  region $\Omega \in  [0,1]^n$ with  smooth
boundary.

For primes $p$ such that $p^n-1$ has a certain prescribed arithmetic
structure M.-C.~Chang~\cite{Chang} obtained nontrivial results which
hold for all primitive roots, rather than on average. Such primes are
rather sparse but one can show that there are infinitely many
of them.

  Now, for an integer $k \ge 1$
and a given region $\Omega \in  [0,1]^n$ with  smooth
boundary, we denote by $N_{\vartheta,k}(\Omega)$ the number of
points
\begin{equation}
\label{eq:Points-k}
\frac{1}{p} \vec{a}_{m^k}, \qquad m =1,\ldots, p^n-1,
\end{equation}
which belong to $\Omega$.

As before, we define  $\rho(k)$ by~\eqref{eq:rho}.

\begin{theorem}
\label{thm:Ergod-k} For any primitive root $\vartheta \in \F_{p^n}$ and any
  region $\Omega \in [0,1]^n$ with piecewise  smooth
boundary, we have
$$
N_{\vartheta,k}(\Omega) = p^n \vol \Omega + O\(  p^{n- \rho(k) + o(1)} \).
$$
\end{theorem}

\begin{proof}
Arguing as in the proof of Theorem~\ref{thm:Ergod-f}, and
using Lemma~\ref{lem:Exp Sum-k} instead of Lemma~\ref{lem:Double Exp Sum},
we easily deduce that the box discrepancy of the set~\eqref{eq:Points-k}
can be estimated as
  $O\( p^{-\rho(k) n + o(1)}\)$. Now applying Lemmas~\ref{lem:Omega-eps}
and~\ref{lem:Omega Discr} we conclude the proof.
\end{proof}

We remark that similar results can be obtained in a more
general situation (for example without  the request that $\vartheta$
is a primitive root). We however follow the settings which exactly
correspond to those of Section~2.A of~\cite{Arn4}.

\subsection{Some other conjectures and open questions}

In Section~2.B of~\cite{Arn4} a conjecture is made which essentially means
that the consecutive values $\vartheta^m$, $\vartheta^{m+1}$
are independently distributed. It is easy to see that this is
incorrect. For example, if $n =1$ and $\vartheta = 2$
is primitive root modulo $p$ (see~\cite{Hool}) then
if $\vartheta^m \in (p/4, p/2)$ then  $\vartheta^m \in (p/2, p)$,
which happens for about $p/4$ values of $m =1, \ldots, p-1$,
while a conjecture given in Section~2.B of~\cite{Arn4} predicts that
this should happen for about $3p/16$ values of $m$.

Sections~2.C and 2.D of~\cite{Arn4} contain a number of interesting
questions about the geometric properties of the set of
points~\eqref{eq:Points}.
We remark that the bound  $O\(M^{-1} p^{n/2} (\log p)^{n+1}\)$ on
the box discrepancy of~\eqref{eq:Points}
obtained in Section~\ref{sec:discr},
immediately implies that any aligned cube $[\alpha, \beta]^n$ inside of the
unit cube $[0,1]^n$ with the side length $\beta - \alpha > C M^{-1/n} p^{1/2}
(\log p)^{1 + 1/n}$, for an appropriate constant $C > 0$,
  contains at least one   point~\eqref{eq:Points}. This immediately
implies upper bounds of the same order on the largest
distance between the points~\eqref{eq:Points} and on the largest radius
of a ball inside of $[\alpha, \beta]^n$ which
does not  contain any points~\eqref{eq:Points}.
Moreover, using some standard modifications, see~\cite{Chalk},
one can drop the logarithmic factor from these bounds.

\section{Dynamical systems in residue rings}

\subsection{Preliminaries}

Given an integer $g\ge 2$ with $\gcd(g, \ell)$,
V.~I.~Arnold~\cite{Arn1,Arn3,Arn5,Arn6}  suggests to consider
the dynamical properties of iterations of
the map $x \mapsto g x$ in the residue ring $\Z_\ell$
(which is equivalent to studying the residues   $g^m \pmod \ell$,
in particular to studying the multiplicative order $t_g(\ell)$).
In particular, in the papers~\cite{Arn1,Arn3,Arn5,Arn6} a number of
suggestions have been made about the average orbit
length of this and several similar dynamical systems.

We remark that indeed if the orbit length is sufficiently large
then, following the standard scheme,
one can derive some analogues of~\eqref{eq:Conj 2A-strong}
from well known bounds of exponential  sums~\cite{KoSh,Kor,LN,Nied1,Nied5}.

Here we provide a brief guide to the literature and
demonstrate that many existing techniques are suitable
for studying these questions and in fact imply that some conjectures
of~\cite{Arn1,Arn3,Arn5,Arn6}, based on numerical calculations,
need some further adjustments.

Throughout this section, the  implied constants
in the Landau symbol `$O$' and in the Vinogradov symbols
`$\ll$' and  `$\gg$'
may occasionally, where  obvious, depend on $g$, and are absolute otherwise
(we recall that $U\ll V$ and $V \gg U$  are both equivalent to the
inequality $U = O(V)$).

\subsection{Analytic number theory background}

Let $\pi_g(x)$ denotes the number of primes $p \le x$, such that
$g$ is a primitive root modulo $p$.

We recall the following celebrated result of Hooley~\cite{Hool}:

\begin{lemma}
\label{lem:Artin} Under the
{\it Extended Riemann Hypothesis\/},   for every  integer $g$ which is
not a perfect square, there exists a constant $A(g) > 0$ such that
$$
\pi_g(x) \sim A(g) \frac{x}{\log x}.
$$
\end{lemma}

Let $\pi(x;k,a)$ denote the number of primes $p \le x$ with
$p \equiv a \pmod k$. We need the following relaxed version of the
{\it  Brun--Titchmarsh\/} theorem,
  see Theorem~3.7 in Chapter~3 of~\cite{HR}.

\begin{lemma}
\label{lem:BT}
For any integers $k,a \ge 1$ with
$1 \le k < x $,  the bound
$$
\pi(x;k,a)  = O\( \frac{ x}{ \varphi(k) \log (3x/k)}\)
$$
holds.
\end{lemma}

Let $\cP$ be the set of prime numbers.
The following estimate can be derived via partial summation
from Lemma~\ref{lem:BT},
see, for example, the proof of Theorem~3.4 in~\cite{EGPS}.

\begin{lemma}
\label{lem:Mertens in A.P.}
For any integer  $f \ge 1$  the bound
$$
  \sum_{\substack{p \in \cP , \, f^2 \le p \le x\\p
\equiv 1
\pmod f}}
\frac{1}{p\log(3x/p)} \ll \frac{\log \log x}{\varphi(f) \log x}
$$
holds.
\end{lemma}

\begin{proof} Let $h = \fl{2\log f}$,  $H = \rf{\log x}$. Then 
\begin{eqnarray*}
\lefteqn{ \sum_{\substack{p \in \cP , \, f^2 \le p \le x\\p
\equiv 1
\pmod f}}
\frac{1}{p\log(3x/p)} \le 
 \sum_{j = h}^H\sum_{\substack{p \in \cP , \,   e^j \le p \le    e^{j+1}\\p
\equiv 1
\pmod f}}
\frac{1}{p\log(3x/p)} }\\
& & \qquad  \ll 
 \sum_{j = h}^H \frac{1}{e^j\log(3x e^{-j-1})}\cdot \frac{e^j}{\varphi(f) j}
\ll  \frac{1}{\varphi(f)} \sum_{j = 1}^H \frac{1}{j\(\log(3x) - j -1\)} .
\end{eqnarray*}
The result now follows.
\end{proof}

\subsection{Average multiplicative order}

Here we show  that
the conjecture~\eqref{eq:Conj Av Period} is wrong and in fact
$T_g(L)$ grows
faster.

\begin{theorem}
\label{thm:AverPer} Under the
{\it Extended Riemann Hypothesis\/},
for every  integer $g$ which is
not a perfect square, there exists a constant $C(g) > 0$ such that
$$
T_g(L)\ge \frac{L}{\log L}
\exp\(C(g)(\log \log \log L)^{3/2}\).
$$
\end{theorem}

\begin{proof} Let $\cP_g$ be the set of  $p \in \cP$
for which  $g$ is a primitive root modulo $p$.
Let us put
$$
Q = \exp(\sqrt{\log L}).
$$
For an integer $k \ge 2$ we  consider the set $\cL_g(k,L)$ of positive integers
$\ell \in [L/2, L]$ of the form
$\ell = p_1\ldots p_k$ where $p_i \in \cP_q$ and $p_i \ge Q$, $i=1, \ldots, k$.

By Lemma~\ref{lem:Artin}, considering only those integers $ \ell \in
\cL_g(k,L)$  for which
$p_1, \ldots, p_{k-1} \le  (L/2Q)^{1/(k-1)}$, and thus $L/2p_1 \ldots
p_{k-1} \ge Q$, we have
\begin{eqnarray*}
\lefteqn{
\# \cL_g(k,L) \ge  \frac{1}{k!}
\sum_{\substack{p_1, \ldots, p_{k-1} \in \cP_q\\
Q \le p_1, \ldots, p_{k-1}  \le  (L/2Q)^{1/(k-1)} }}
\, \sum_{\substack{p_k \in
\cP_q\\
  L/2p_1 \ldots p_{k-1}  \le p_k \le L/p_1 \ldots p_{k-1} }} 1}\\
& &=   \frac{(A(g) + o(1))L}{2k!} \sum_{\substack{p_1, \ldots,
p_{k-1} \in \cP_q\\
Q \le p_1, \ldots, p_{k-1}  \le  (L/2Q)^{1/(k-1)}}} \frac{1}{p_1 \ldots p_{k-1}
\log(3L/p_1 \ldots p_{k-1})}\\
& &\ge
\frac{(A(g) + o(1)) L}{2k! \log L}  \sum_{\substack{p_1, \ldots, p_{k-1} \in
\cP_q\\  Q \le p_1, \ldots, p_{k-1}  \le  (L/2Q)^{1/(k-1)}}}
\frac{1}{p_1 \ldots
p_{k-1}}\\ & & =  \frac{(A(g) + o(1))L}{2k!}
\(\sum_{\substack{p \in \cP_q\\ Q \le  p \le (L/2Q)^{1/(k-1)}}}
\frac{1}{p}\)^{k-1}.
\end{eqnarray*}
By partial summation, we derive from  Lemma~\ref{lem:Artin} that
$$
  \sum_{\substack{p \in \cP_q\\  Q \le  p \le (L/2Q)^{1/(k-1)} }}
\frac{1}{p} \sim
A(g) (\log \log (L/2Q)^{1/(k-1)}  - \log \log Q) \sim  \frac{A(g)}{2}
\log \log L,
$$
uniformly for $k$ with $\log k = o(\log \log L)$.
Therefore, uniformly for $k$ with $\log k = o(\log \log L)$
$$
\# \cL_g(k,L) \ge \frac{ \( A(g)    + o(1)\)^{k} L (\log \log
L)^k}{2^{k} k!\log L}.
$$

We now denote by $\cQ_d(k,L)$ the set of  positive integers $\ell\le
L$ of the form
$\ell = p_1\ldots p_k$ where $p_i \in \cP$ are distinct primes,
$p_i \ge Q$, $i=1, \ldots, k$, and
$$
\max_{1 \le i < j \le k} \gcd(p_i-1,p_j-1) = d.
$$

For $d \le Q^{1/2}$ we use  Lemma~\ref{lem:BT} to derive
\begin{eqnarray*}
\# \cQ_d(k,L) &\le & \frac{1}{(k-2)!}
\sum_{\substack{p_1, \ldots, p_{k-2} \in \cP \\   p_1, \ldots, p_{k-2}
\ge Q}} \sum_{\substack{p \in \cP \\    p \ge Q \\ p \equiv 1 \pmod d}}
\sum_{\substack{q \in \cP \\   Q \le q \le L/pp_1\ldots p_{k-2} \\ q
\equiv 1 \pmod d}} 1 \\
& \ll & \frac{L}{
(k-2)! \varphi(d)}
\sum_{\substack{p_1, \ldots, p_{k-2} \in \cP \\   p_1, \ldots, p_{k-2}
\ge Q}} \\
& &\qquad \qquad \sum_{\substack{p \in \cP \\  Q \le p \le
L/p_1\ldots p_{k-2} \\
p \equiv 1
\pmod d}}
\frac{1}{pp_1\ldots p_{k-2} \log(3L/pp_1\ldots p_{k-2}d)} .
\end{eqnarray*}
Applying Lemma~\ref{lem:Mertens in A.P.} with $f=d$ (once) and theh $f=1$
($k-2$ times), we
obtain, that for
$d
\le Q^{1/2}$
$$
\# \cQ_d(k,L)  \ll  \frac{L (\log \log L)^{k-1}}{(k-2)! \varphi(d)^2\log L}.
$$
Since $\varphi(d) \gg d/\log \log d$,  see Theorem~328
of~\cite{HardyWright}, we conclude that for any $D > 0$,
$$
\sum_{Q^{1/2} \ge d > D} \# \cQ_d(k,L) \le D^{-1 + o(1)} \frac{L
(\log \log L)^{k-1}}{(k-2)!
\log L}.
$$
Also using the trivial bound
$$
\# \cQ_d(k,L)  \ll \sum_{\substack{p_1,p_2 \in \cP, \\ p_1,p_2
\le L \\p_1 \equiv p_2
\equiv 1 \pmod d}} \frac{L}{p_1p_2} \le
L \(\sum_{\substack{n
\le L \\n
\equiv 1 \pmod d}} \frac{1}{n} \)^2 \ll \frac{L (\log L)^2 }{d^2}.
$$
for $d > Q^{1/2}$ we derive that
$$
\sum_{d > Q^{1/2}} \# \cQ_d(k,L) \le Q^{-1/2} L  (\log L)^2.
$$
Hence,  for $D \le Q^{1/2}$
$$
\sum_{d > D} \# \cQ_d(k,L)  \le   D^{-1 + o(1)}
\frac{  L(\log \log L)^{k-1}}{(k-2)!
\log L}.
$$
Finally, let $\cS(L)$ be the set of  positive integers $\ell\le L$
such that $p^2 | \ell$ for some $p\ge Q$.
Clearly
$$
\# \cS(L) \le \sum_{\substack{p \in \cP \\   p \ge Q}} \fl{L/p^2} \ll L/Q.
$$

Therefore, for $D_k = 3^k$, for  the set
$$
\cR_g(k,L) = \cL_g(k,L) \setminus \(\bigcup_{d > D_k}  \cQ_d(k,L) \cup \cS(L)\)
$$
we have
$$
\# \cR_g(k,L) \ge  \frac{ \( A(g)    + o(1)\)^{k} L (\log \log
L)^k}{2^{k} k!\log L}, 
$$
provided  $\log k = o(\log \log L)$.

On the other hand, for every $\ell = p_1\ldots p_k \in \cR_g(L)$ we have
\begin{eqnarray*}
t_g(\ell) & = & \mathrm{lcm}(t_g(p_1), \ldots, t_g(p_k)) =
\mathrm{lcm}(p_1-1,\ldots, p_k-1) \\
& \ge & \prod_{i=1}^k {(p_i-1)} \prod_{j =1}^{i-1} \frac{1}{\gcd
(p_i-1, p_j-1)}\\
& \ge & \frac{(p_1-1)\ldots (p_k-1)}{D_k^{k^2/2}}\ge \frac{p_1\ldots p_k }
{2^kD_k^{k^2/2}} \ge
\frac{L}{2^{k+1}D_k^{k^2/2}}\\
& \ge &
\frac{L}{2^{k+1}3^{k^3/2}} \gg \frac{L}{3^{k^3}}.
\end{eqnarray*}
Therefore,  under the above condition  on  $k$, we derive
\begin{eqnarray*}
T_g(L) & \ge &
\frac{1}{L} \sum_{\substack{\ell \in \cR_g(k,L) \\ \gcd(g,\ell)=1}}^L t_g(\ell)
\ge   \frac{ \# \cR_g(k,L) }{3^{k^3}}
  \ge  \frac{ \( A(g)    + o(1)\)^{k} L (\log \log L)^k} { k! 3^{k^3} \log
L}\\ & =& \frac{  L (\log \log L)^k \exp(O(k^3))} {\log L}.
\end{eqnarray*}
Taking
$$
k =  \fl{c(g) \sqrt{\log \log \log L}}
$$
for an appropriate constant $c(g)>0$, depending only on $g$
(which guarantees that the term $(\log \log L)^k$ exceeds,
say,  the square of
factor $\exp(O(k^3))$)
we  finish the proof.
\end{proof}

As we have remarked, we believe that the bound of Theorem~\ref{thm:AverPer}
is not tight and in fact a stronger bound~\eqref{eq:AverPerTight}
can be derived by using the method of~\cite{MaPom}, modified in
a similar way as that of~\cite{BFLPS} to deal only with special primes,
see also~\cite{Ford1}.

It is a very interesting question to obtain more precise
information about the behaviour of $T_g(L)$, for example
to establish whether~\eqref{eq:T and lambda} is correct. It is possible
that the method of~\cite{EPS} combined with the methods and results
of~\cite{KurPom,LiPom}  are  able to handle this task. In fact, even already
existing results of~\cite{EPS,KurPom,LiPom}, without any modifications or
adjustments,  may shed light on  many issues  risen by
V.~I.~Arnold in~\cite{Arn1,Arn3,Arn5,Arn6}.

We also remark that a dual question about the the average value
$$
\widetilde{T}(\ell) = \frac{1}{\varphi(\ell)} \sum_{\substack{g =1\\
\gcd(g,\ell)=1}}^\ell t_g(\ell)
$$
is studied  in~\cite{GKLLS,Luc2,LucShp}. In particular,
in~\cite{LucShp}  one can also find
various  upper and lower bounds on  $\widetilde{T}(\ell)$, while its
behaviour on
special sequences is considered in~\cite{GKLLS,Luc2}.

It is well known
that if an integer $g>1$ is fixed then for  any function $\varepsilon(x)$
with $\varepsilon(x) \to 0$ as $x \to \infty$,
for almost all primes $p$ the bound
$ t_g(p) \ge p^{1/2 + \varepsilon(p)}$ holds,
see~\cite{ErdMur,Ford2,IndlTim,Papp} for various improvements of this result.
For almost all integers $\ell$, similar type bounds are given in~\cite{KurPom}.

It is clear that when $g$ varies, $t_g(\ell)$ runs through divisors of
$\lambda(\ell)$. In fact through all the divisors of $\lambda(\ell)$.
Accordingly the question about the behaviour of $\tau(\lambda(\ell))$
becomes of interest (where $\tau(k)$ is the number of all
integer positive divisors of $k \ge 1$). Partially motivated by this
relations, F.~Luca and C.~Pomerance~\cite{LucPom} obtained tight
bounds on the average
value of $\tau(\lambda(\ell))$.

\subsection{Average additive order}

V.~I.~Arnold~\cite{ Arn5} also asks about the average
period $Q(\ell)$ of the map $x \mapsto x+a$ in $\Z_\ell$ taken over
all  $a \in \Z_\ell$. This function, which can be expressed as
the following sum over the divisors of $\ell$
$$
Q(\ell) = \frac{1}{\ell} \sum_{d|\ell} d \varphi(d),
$$
has been studied in detail in~\cite{GKLLS,Luc2}.

In particular, it is shown in Theorem~3.1 of~\cite{GKLLS}
that
$$
\frac{1}{L} \sum_{\ell \le L} Q(\ell) = \frac{3\zeta(3)}{\pi^2} L
+ O(  \(\log L)^{2/3} (\log \log L)^{4/3}\)
$$
where $\zeta(s)$ is the {\it Riemann $\zeta$-function\/}.
This gives a more precise and explicit form of the
assertion made in~\cite{Arn5} that on average $Q(\ell)$ grows linearly.
Clearly, for any prime $p$, we have $\widetilde{T}(p) = Q(p-1)$, thus
it is also interesting to
study  $Q(\ell)$ on this and some special sequences of $\ell$,
see~\cite{GKLLS}.
More results on arithmetic properties of $Q(\ell)$ have been established
by F.~Luca~\cite{Luc2}.

\subsection{Average divisor}

We note that well known bounds of number theoretic functions
implies that the assertion  made
in~\cite{Arn5}   that the ``average divisor'' of $\varphi(\ell)$ is
$$
d(\varphi(\ell)) = \frac{\sigma(\varphi(\ell))}{\tau(\varphi(\ell))}
\sim \frac{\ell^\beta} {\log \ell}
$$
(where $\sigma(k)$ is the sum of all
integer positive divisors of $k \ge 1$)
with some $\beta = 0.96 \pm 0.02$ is false. As it follows from
the classical number theoretic bounds
$$
2 \le \tau(k) \le 2^{(1+o(1)) \log k/ \log \log k},   \qquad
k + 1\le \sigma(k) =
O(k \log \log k),
$$
see Theorems~317 and~323 of~\cite{HardyWright}, respectively,
for any sufficiently large $k$ we have
$$
d(k) = \frac{\sigma(k)}{\tau(k)}
\sim k^{1 + o(1)}.
$$
Finally,  we mention that the suggestion made in~\cite{Arn5} that
the average value of $d(k)$ behaves like
$$
D(K) = \frac{1}{K}\sum_{k=1}^K d(k) \sim \frac{3K}{2 \log K}
$$
is wrong too. Since $d(k)$ is a multiplicative function, so is
$d(k)/k$, which also satisfies the conditions of the
{\it Wirsing theorem\/}, see~\cite{Wirs}. Thus one can easily
show that in fact
$$
D(K)    \sim \kappa \frac{K}{(\log K)^{1/2}}
$$
for some absolute constant $\kappa = 0.4067\ldots$, see~\cite{BEPS} for
this and some other results on the properties of the average
divisor, including an asymptotic expansion of $D(K)$. However,
  the
question on the average value of $d(\varphi(\ell))$ is
harder and is of ultimate interest.

It can also be relevant to mention that  average divisor $d(k)$
takes integer values  for almost
all integers $k\ge 1$  but is almost never a divisor of $k$,
see~\cite{BEPS,Luc1,Spiro} and the references therein.
Certainly similar questions about $d(\varphi(\ell))$ are very
natural and interesting.

\section{Repeated squaring and other nonlinear\break transformations}

Using~\cite{BFGS,Bour2,FHS1,FHS2,FrKoSh,FrSh} one can also easily derive
various uniformity
of distribution results for the vectors  $\vec{a}_{e^m}$
where $e\ge 2$ is a fixed integer.  Alternatively,  these results can be
interpreted as results about orbits of repeated powering
$x \mapsto x^e$. In particular, with $e
=2$ one can study the distribution of elements in orbits of repeated squaring
$x \mapsto x^2$
in finite fields and rings, see~\cite{Arn2,Arn3} where the
corresponding dynamical
system is outlined. The results of~\cite{FHS1,FHS2,FrKoSh,FrSh}
show that if the orbit is long enough then the vectors  $\vec{a}_{e^m}$
are uniformly distributed. Using the bound of J.~Bourgain~\cite{Bour2},
one can
consider very short orbits (and not necessary fixed values of $e$)
although the bounds obtained within this approach are less explicit.

These results are complemented by the estimates on
the orbit lengths  of such transformations which are obtained
in~\cite{FPS,MartPom} and which show that these
orbit lengths   tend to be
large (and close to their largest possible values).

The distributional properties of
dynamical systems generated by general non-linear transformations
$x \mapsto f(x)$ where $f$ is a rational function,
over a finite field  or a residue ring,  have been extensively studied in the
literature as well,
see~\cite{ElShWi,ElWi,Nied6,NiSh5,NiWint3}
and the references therein.
As in the case of repeated powering all these results indicate that
if the orbit is long enough then its elements are uniformly distributed.
On the other hand, these results are still missing their essential
counterpart, namely estimates on the orbit length.  Obtaining such
estimates (for general or specific functions $f$) is a very important
open question.

\section{Further remarks and extensions}

Results of a different flavour
but also describing the distribution of powers of primitive elements
 modulo a prime $p$ are given in~\cite{RudZah}.

We have already remarked that analogues of our results hold
for an arbitrary $\vartheta$, not necessarily a primitive root,
provided the multiplicative order of $\vartheta$ is large enough.
For more general formulations
of  Lemma~\ref{lem:Double Exp Sum}  and  Lemma~\ref{lem:Exp Sum-k},
see~\cite{Gar1,GarKar} and~\cite{FHS1},
respectively.

Analogues of the
bound~\eqref{eq:Classic Exp Sum} and  Lemmas~\ref{lem:Double Exp Sum}
and~\ref{lem:Exp Sum-k} are also known for residue rings $\Z_\ell$,
see~\cite{KoSh,Kor,LN}
and~\cite{FHS2,FrKoSh},  respectively.
Thus one can study orbits of $\vartheta^m$ in
residue rings as well. In particular, it is well known that one can
use these bounds to
obtain various uniformity of distribution results suggested in
Section~2 of~\cite{Arn3}.

Furthermore, the dynamical system corresponding to the repeated
squaring of a  unimodular matrix  has been considered
in~\cite{Arn2}. Using the results and methods of~\cite{GoGuSh},
one can prove various uniformity of distribution properties
of orbits of such dynamical systems. Accordingly,
the results of~\cite{BFGS,HeSh,LanShp} can be used to obtain similar
statements for analogues of the above dynamical systems on
elliptic curves over finite fields.

\section*{Acknowledgements}

  The author wishes to thank Florian Luca, Harald
Niederreiter, Carl Pomerance, John Roberts  and Franco Vivaldi for  useful
discussions, additional references and encouragement  to write this paper.
This work was supported in part by ARC grant  DP0556431.

\end{document}